\documentclass[a4paper,12pt]{article}

\usepackage{amsmath}
\usepackage{amssymb}
\usepackage[utf8]{inputenc}
\usepackage[english]{babel}
\usepackage{graphicx} % for pdf, bitmapped graphics files
\usepackage{color}
\usepackage{subfigure}
\usepackage{enumitem}
\usepackage{tikz}
\usepackage{url}

\usepackage{bm}

\title{On spectral properties for graph matching and graph isomorphism problems}
\newcommand{\Pe}{\mathcal{P}}

\newcommand{\ones}{\mathbf{1}}

\newcommand{\bb}[1]{\bm{\mathrm{#1}}}

\newcommand{\ds}{\mathcal{D}}

\def\reals{\ensuremath{\mathbb{R}}}

\newtheorem{theorem}{Theorem}
\newtheorem{lemma}[theorem]{Lemma}
\newtheorem{corollary}[theorem]{Corollary}
\newtheorem{proposition}[theorem]{Proposition}
\newtheorem{conjecture}[theorem]{Conjecture}

\begin{document}
\author{Marcelo Fiori \\ Universidad de la Rep\'ublica, Uruguay\\ mfiori@fing.edu.uy \\ \\ Guillermo Sapiro \\ Duke University\\ guillermo.sapiro@duke.edu}
\date{}
%\shortauthorlist{Marcelo Fiori and Guillermo Sapiro} %%% for verso running head

\maketitle

\begin{abstract}
{Problems related to graph matching and isomorphisms are very important both from a theoretical and practical perspective, with applications ranging from image and video analysis to biological and biomedical problems. The graph matching problem is  challenging from a computational point of view, and therefore different relaxations are commonly used. Although common relaxations techniques tend to work well for matching perfectly isomorphic graphs, it is not yet fully understood under which conditions the relaxed problem is guaranteed to obtain the correct answer.

In this paper we prove that the graph matching problem and its most common convex relaxation, where the matching domain of permutation matrices is substituted with its convex hull of doubly-stochastic matrices, are equivalent for a certain class of graphs, such equivalence being based on spectral properties of the corresponding adjacency matrices. We also derive results about the automorphism group of a graph, and provide fundamental spectral properties of the adjacency matrix.}

%\red{abstract}
\end{abstract}

\vspace*{1cm}

\section{Introduction}

The theoretical and computational aspects behind graph isomorphisms and graph matching have been a great challenge for the scientific community for a long time. Maybe the easiest problem to state from this category is the graph isomorphism problem, which consists in determining whether two given graphs are isomorphic or not, that is, if there exists a bijection between the vertex sets of the graphs, preserving the edge structure. Besides the theoretical analysis, the graph isomorphism problem is also very interesting from the computational complexity point of view, since its complexity class is still unsolved: it is one of the few problems in NP not yet classified as P nor NP-complete \cite{ConteReview}. 

The concept of graph automorphism, and its related properties, is closely connected to the graph isomorphism problem. An automorphism of a graph is a mapping from its vertex set onto itself, preserving the connectivity structure. The set of automorphisms forms a group under the composition operation. Of course, the identity map is always an automorphism, and when this is the only element in the group, we say that the graph has a trivial automorphism group. From the computational complexity point of view, computing the automorphism group is at least as difficult as solving the graph isomorphism problem. 

The last problem we wish to discuss here is the so-called graph matching problem, which consists in finding an isomorphism between two graphs, and it is therefore harder than the graph isomorphism problem. Specifically, let $G_A$ and $G_B$ be two graphs with $n$ vertices, and let $\bb{A}$ and $\bb{B}$ be their corresponding adjacency matrices. A common statement of the graph matching problem is to find the correspondence between the nodes of $G_A$ and $G_B$ which minimizes some matching error. 
In terms of the corresponding adjacency matrices $\bb{A}$ and $\bb{B}$, which encode the graph connectivity,  this corresponds to finding a matrix $\bb{P}$ in the set of permutation matrices $\Pe$, such that it minimizes a given distance between $\bb{A}$ and $\bb{PBP^T}$. A common choice is the Frobenius norm $||\bb{A} - \bb{PBP^T}||_F^2$, and then the graph matching problem can be formally stated as

\begin{equation}
\label{eq:GM}
\tag{$P_1$}
\min_{\bb{P}\in\Pe}||\bb{A} - \bb{PBP^T}||_F^2 = \min_{\bb{P}\in\Pe}||\bb{AP} - \bb{PB}||_F^2.
\end{equation}

Although polynomial algorithms have been developed for a few special types of graphs, like trees or planar graphs for example \cite{ConteReview}, the combinatorial nature of the permutation search makes this problem NP in general. As such, there are several and diverse techniques addressing the graph matching problem, including spectral methods \cite{umeyama} and relaxations techniques \cite{fiori2013nips,FAQ,Zaslavskiy2009}.

%In this paper we focus on a particular relaxation technique, which consists in relaxing the feasible set (the set of permutation matrices) to its convex hull. By virtue of the Birkhoff-von Neuman theorem, the convex hull of $\Pe$ is the set of doubly stochastic matrices $\ds$, which is the set of matrices with non-negative entries such that each row and column sum up one: $\ds=\{\bb{M}\in \reals^{p\times p} : \bb{M}_{ij}\geq 0, \bb{M}\ones=\ones, \bb{M}^T\ones=\ones\}$.

In this paper we focus on a particular and very common relaxation technique, which consists in relaxing the feasible set (the set of permutation matrices) to its convex hull. By virtue of the Birkhoff-von Neuman theorem, the convex hull of $\Pe$ is the set of doubly stochastic matrices $\ds=\{\bb{M}\in \reals^{n\times n} : \bb{M}_{ij}\geq 0, \bb{M}\ones=\ones, \bb{M}^T\ones=\ones\}$, that is, the set of matrices with non-negative entries such that each row and column sum up to one.

%, $\uno$ being the $p$-dimensional vector of ones. %Observing that for permutation matrices $||\bb{A} - \bb{PBP^T}||_F^2 = ||\bb{AP} - \bb{PB}||_F^2$, 
The relaxed version of the  problem is then
\begin{equation}
\label{eq:RGM}
\tag{$P_2$}
\hat{\bb{P}}=\arg\min_{\bb{P}\in\ds}||\bb{AP} - \bb{PB}||_F^2,
\end{equation}
which is a convex problem. However, the resulting $\hat{\bb{P}}$ is a doubly stochastic matrix and not necessarily a permutation matrix, or in general the solution to \eqref{eq:GM}. %The final node correspondence is obtained as the closest permutation matrix to $\hat{\bb{P}}$: $\bb{P}^* = \arg\min_{\bb{P}\in\Pe}||\bb{P}-\hat{\bb{P}}||_F^2$, which is a linear assignment problem that can be solved in $O(p^3)$ by the Hungarian algorithm \cite{hungarian}. However, this last step lacks any guarantee about the graph matching problem itself. This approach will be referred to as QCP for \textit{quadratic convex problem}.

Indeed, since the feasible set of problem \eqref{eq:RGM} is the convex hull of the feasible set of problem \eqref{eq:GM}, every solution of the first problem is also a solution of the relaxed graph matching problem.  A very important question is under which hypothesis the solution set of these two problems \eqref{eq:GM} and \eqref{eq:RGM} coincide. It is easy to see that, if there are two permutation matrices that solve the graph matching problem, then every matrix on the straight line joining them is a solution of problem \eqref{eq:RGM}, since this problem is convex. Therefore, the least that one should ask for these two problems to be equivalent is for the solution of \eqref{eq:GM} to be unique. When the two graphs are isomorphic, this is equivalent to asking for the automorphism group of the graphs to be the trivial group.

A probabilistic analysis of this equivalence between the original and the relaxed graph matching problems is provided in \cite{vince}. %, where a non-convex relaxation is also studied with very encouraging theoretical results.
 Basically, the authors prove that, when two graph are correlated (but not necessarily isomorphic), then the unique solution of a non-convex relaxation is almost always the correct permutation matrix; while on the other hand, the underlying alignment is almost always not a solution of the commonly used convex relaxation, where the permutation set is replaced by the doubly-stochastic set as above. 

On the other hand, in \cite{alex} the authors prove the equivalence of the original graph matching problem and a relaxed version for a particular kind of graphs which they call \textit{friendly}, based on spectral properties. In this work, we extend these results, proving the (deterministic) equivalence for a larger set of graphs, and also shedding light on some new spectral graph properties.

%\newpage

\section{Main result}
\label{sec:thms}

In this section, we consider two isomorphic graphs $G_A$ and $G_B$ with $n$ vertices each, and adjacency matrices $\bb{A}$ and $\bb{B}$ respectively. Let $\bb{P_o} \in \Pe$ be the permutation matrix associated to the isomorphism between the two graphs, that is, $\bb{B} = \bb{P_oAP_o^T}$.

Since the graphs considered here are isomorphic, then the minimum (either over $\ds$ or $\Pe$) of $\|\bb{AP}-\bb{PB}\|_F^2$ is zero, and it is achieved (at least) at $\bb{P_o}$. Both problems can be then re-stated as solving the set of linear equations $\bb{AP}=\bb{PB}$ over $\bb{P}\in \Pe$ or $\bb{P} \in \ds$.

  Now, consider that by the simple change of variables $\bb{Q}=\bb{PP_o}$. Then for any solution $\bb{P}$ to the relaxed problem \eqref{eq:RGM}, it holds that
  \begin{equation}
  \bb{AP}=\bb{PB} \Longleftrightarrow \bb{AP} = \bb{PP_oAP_o^T} \Longleftrightarrow \bb{APP_o} = \bb{PP_oA} \Longleftrightarrow \bb{AQ}=\bb{QA}.
  \end{equation}

Note that the change of variables is a multiplication by a permutation matrix, and hence the set of doubly stochastic matrices is invariant under this mapping. Therefore, any solution to $\bb{AQ}=\bb{QA}$ over $\bb{Q}\in \ds$ leads, via the change of variables, to a solution of $\bb{AP}=\bb{PB}$ with $\bb{P}\in \ds$. This allows us to state the equivalency between both problems \eqref{eq:GM} and \eqref{eq:RGM} using only one of the adjacency matrices. Specifically, the problem $\bb{AQ}=\bb{QA}$ with $\bb{Q} \in \ds$ has a trivial solution $\bb{Q}=\bb{I}$, which corresponds to the solution $\bb{P_o}$ of the problem $\bb{AP}=\bb{PB}$ with $\bb{P} \in \ds$. Then the matrix $\bb{P_o}$ will be the unique solution of problem \eqref{eq:RGM} if and only if the identity is the unique solution of $\bb{AQ}=\bb{QA}$ with $\bb{Q} \in \ds$.

Now, since $\bb{A}$ is a symmetric matrix, we can consider its spectral decomposition $\bb{A}=\bb{UDU^T}$, where $\bb{D}$ is a diagonal matrix containing the eigenvalues and $\bb{U}$ is an orthonormal matrix containing the eigenvectors as columns, denoted as $u_i$, for $i=1\ldots n$.

\begin{figure}[h!]
\centering
\includegraphics[width=0.49\textwidth]{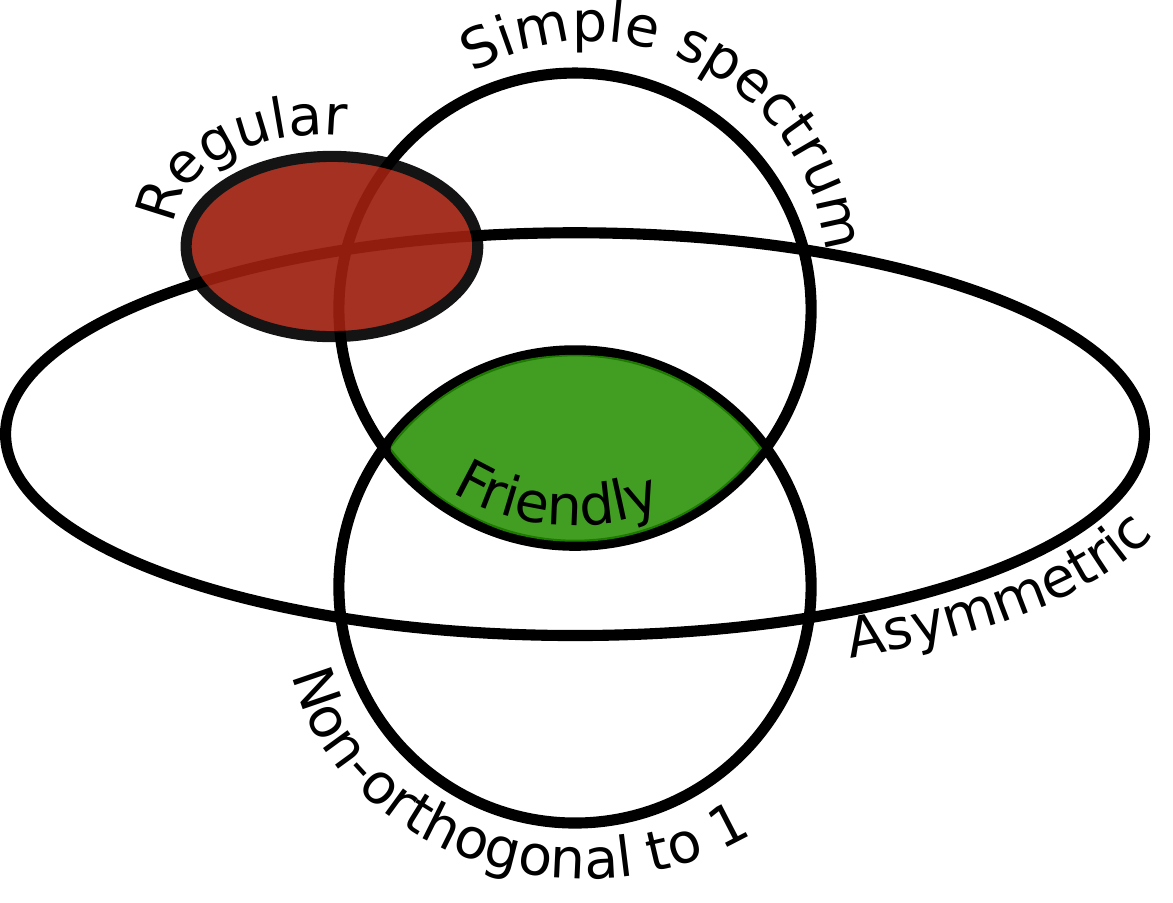} 
\caption{Graph classes and equivalence of problems \eqref{eq:GM} and \eqref{eq:RGM}. Graphs with trivial authomorphism group are represented in the \textit{asymmetric} set (we know that graphs for which $(P_1)$ and $(P_2)$ are equivalent are inside this set), graphs whose adjacency matrices have no repeated eigenvalues are represented as the \textit{simple spectrum} set, \textit{non-orthogonal to $\bb{1}$} means that no eigenvector $u_i$ verifies $u_t^T\ones=0$, and the \textit{regular} circle contains regular graphs. Graphs in the intersection of \textit{simple spectrum} and \textit{non-orthogonal to $\bb{1}$} are called \textit{friendly} graphs, and here the equivalence of problems \eqref{eq:GM} and \eqref{eq:RGM} holds \cite{alex}. These problems are not equivalent for regular graphs. A key question addressed in this work is how far we can extend the green zone of equivalence inside the asymmetric set.}
\label{fig:bubble}
\end{figure}

The main result of \cite{alex} states that if $\bb{A}$ has no repeated eigenvalues, and no eigenvector $u_i$ is perpendicular to the vector of ones $\ones$, % = (1,1,\ldots,1)$, %i.e., $u_i^T\ones \neq 0$ for all $i=1\ldots n$, 
then problems \eqref{eq:GM} and \eqref{eq:RGM} are equivalent. This is illustrated in Figure \ref{fig:bubble}, where some graph properties are represented. Here \textit{asymmetric} means that the authomorphism group of the graph is trivial, \textit{simple spectrum} means that the adjacency matrix has no repeated eigenvalues, \textit{non-orthogonal to $\bb{1}$} means that no eigenvector $u_i$ verifies $u_t^T\ones=0$, and the \textit{regular} circle contains regular graphs, i.e., graphs such that each vertex has the same number of neighbors. The intersection of \textit{simple spectrum} and \textit{non-orthogonal to $\bb{1}$} graphs is what the authors of \cite{alex} call \textit{friendly} graphs, and they prove the equivalence of problems \eqref{eq:GM} and \eqref{eq:RGM} for this class.

As observed above, a necessary condition for problems \eqref{eq:GM} and \eqref{eq:RGM} to be equivalent, is for the automorphism group of the graph to be the trivial group. However, this condition is not sufficient. Take for instance a regular graph, and denote by $\bb{J}$ the barycenter of the set of doubly stochastic matrices, $\bb{J}=\dfrac{1}{n}\ones\ones^T$. Hence, it is very easy to see that, if $\bb{A}$ is the adjacency matrix of a regular graph, then $\bb{AJ}=\bb{JA}$. Therefore, there is a solution to problem \eqref{eq:RGM} which is not a permutation matrix. Since there are regular graphs with trivial automorphism group (like the Frucht graph \cite{frucht} for instance), then this condition cannot be sufficient. In Figure \ref{fig:bubble}, this is represented with the small red circle, which intersects both the asymmetric and simple spectrum sets (see Section \ref{graph_examples} for examples of graphs in each intersection). In summary, we can hope for problems \eqref{eq:GM} and \eqref{eq:RGM} to be equivalent inside the \textit{asymmetric} set minus the regular graphs. So far, we know from \cite{alex} that this is true for \textit{friendly} graphs, being this until now the largest known class for which the relaxation is equivalent to the original problem.

The next theorems extend the set where these problems are equivalent to a larger set of graphs. Theorem \ref{thik} is stronger than Theorem \ref{th2k}, but we include both proofs for the sake of clarity, since both have pedagogic value.

%With the notation defined above, we have the following
\begin{theorem}
%\begin{thm}
\label{th2k}
If $\bb{A}$ has no repeated eigenvalues (simple spectrum), and there are $k$ eigenvectors $u_i$ such that $u_i^T\ones =0$, each one of these vectors having at least $2k+1$ nonzero entries, 
then problems \eqref{eq:GM} and \eqref{eq:RGM} are equivalent.
%\end{thm}
\end{theorem}

%\noindent \textbf{Proof of Theorem \ref{th2k}}
\noindent \textit{Proof:}
%\begin{proof}

We want to prove that the identity is the unique solution to the problem $\bb{AQ}=\bb{QA}$ for $\bb{Q}\in\ds$. Let us write the equality $\bb{AQ}=\bb{QA}$ in terms of the eigenvector decomposition of $\bb{A}$:

\begin{equation*}
\bb{AQ}=\bb{QA} \Leftrightarrow \bb{UDU^TQ}=\bb{QUDU^T} \Leftrightarrow \bb{U^TUDU^TQU}=\bb{U^TQUDU^TU} \Leftrightarrow \bb{DU^TQU}=\bb{U^TQUD}.
\end{equation*}

Now, let us denote by $\bb{F}$ the new unknown matrix $\bb{F}=\bb{U^TQU}$. The problem can be now stated as
\begin{equation}
\bb{DF}=\bb{FD} \quad , \quad \bb{UFU^T} \in \ds,
\label{problemF}
\tag{$P_F$}
\end{equation}
\noindent and we now want to prove that $\bb{F}=\bb{I}$ is the unique solution of this last problem.

It is easy to see that, since $\bb{D}$ is diagonal with no repeated entries in the diagonal, then $\bb{F}$ has to be diagonal as well in order to commute with $\bb{D}$.

Let us write the conditions for $\bb{UFU^T}$ to be in $\ds$:

%\begin{itemize}[noitemsep]
\begin{itemize}
\item[c1)] $\bb{UFU^T} \ones = \ones$,
\item[c2)] $\bb{UF^TU^T} \ones = \ones$,
\item[c3)] $\left( \bb{UFU^T} \right)_{i,j} \geq 0 \, , \, \forall \, i,j$.
\end{itemize}

Since $\bb{F}$ is diagonal, and in particular $\bb{F}=\bb{F^T}$, then the first two conditions are redundant, and one of them can be eliminated. 
Left-multiplying the first condition by $\bb{U^T}$, we obtain $\bb{FU^T}\ones = \bb{U^T}\ones$, and calling $\mathbf{v}=\bb{U^T}\ones$, condition c1) can be written as $\bb{F}\mathbf{v}=\mathbf{v}$.

Without loss of generality, we can assume that the $k$ eigenvectors $u_i$ satisfying $u_i^T\ones=0$ are the first $k$ columns in $\bb{U}$. Therefore, $\mathbf{v}_i=0$ for $i=1\ldots k$, and $\mathbf{v}_i\neq 0$ for $i=k+1\ldots n$. As $\bb{F}$ is diagonal, the equations  $\bb{F}\mathbf{v}=\mathbf{v}$ can be easily written as $\bb{F}_{i,i}\mathbf{v}_i = \mathbf{v}_i$. When $\mathbf{v}_i\neq 0$, the only way for this equation to hold is when $\bb{F}_{i,i}=1$. This means that $\bb{F}_{i,i}=1$ for $i=k+1\ldots n$, and this is sufficient to guarantee that the first two conditions hold.

For analyzing the third condition, let us decompose the matrix product using that $\bb{F}$ is a diagonal matrix with $\bb{F}_{i,i}=1$ for $i=k+1\ldots n$:

$$\bb{UFU^T} = \sum_{i=1}^{n} u_i\bb{F}_{i,i}u_i^T = \sum_{i=1}^{k} u_i\bb{F}_{i,i}u_i^T + \sum_{i=k+1}^{n} u_iu_i^T .$$

We can now add and subtract $\sum_{i=1}^{k}u_iu_i^T$, leading to 

$$\bb{UFU^T} = \sum_{i=1}^{k} u_i(\bb{F}_{i,i}-1)u_i^T + \sum_{i=1}^{n} u_iu_i^T = I+ \sum_{i=1}^{k} u_i(\bb{F}_{i,i}-1)u_i^T.$$

%\red{Observe that the values $\bb{F}_{i,i}$ must satisfy $0\leq \bb{F}_{i,i}\leq 1$ in order for the entries of $\bb{UFU}^T$ to be between $0$ and $1$.}
Let us denote by $\bb{L} =  \sum_{i=1}^{k} (1-\bb{F}_{i,i})u_iu_i^T$, and therefore $\bb{UFU^T} = \bb{I}-\bb{L}$. 

Observe that the matrix $\bb{L}$ satisfies $\bb{L}\ones=\mathbf{0}$, since every vector $u_i$ participating in the sum satisfies $u_i^T\ones=0$; and note also that all the elements in the diagonal are $\bb{L}_{j,j} \geq 0$, otherwise the corresponding entry of $\bb{UFU^T}$ would be $(\bb{UFU}^T)_{j,j}>1$, violating the doubly stochastic condition. %because they are the result of the sum of positive numbers: $\bb{L}_{j,j} = \sum_{i=1}^{k} (1-\bb{F}_{i,i})(u_{i})_j^2$.

Now, let us assume that there is a solution $\bb{F}$ to problem \eqref{problemF} different from the identity, and let us analyze the corresponding $\bb{L}$ matrix trying to find a contradiction. The condition c3) dictates that $(\bb{I}-\bb{L})_{i,j} \geq 0$ for all $i,j$, therefore the $\bb{L}$ matrix has no positive elements off the diagonal. On the other hand, since $\bb{F}$ is diagonal and we have assumed $\bb{F}\neq \bb{I}$, then at least one of the values $\bb{F}_{i,i}$  ($i \leq k$) is different from $1$. The corresponding eigenvector $u_i$, which has at least $2k+1$ non-zero elements by hypothesis, will be actually used in the summation constructing $\bb{L}$, and therefore this guarantees that at least $2k+1$ elements in the diagonal of $\bb{L}$ are strictly positive.

 Considering then the following just described properties for the $\bb{L}$ matrix:
 
 \begin{itemize}
 \item $\bb{L}\ones = \mathbf{0}$,
 \item $\bb{L}_{i,j} \leq 0$ for all $i\neq j$,
 \item $\bb{L}_{i,i} \geq 0$ for all $i=1\ldots n$,
 \item $\bb{L} = \bb{L}^T$.
 \end{itemize}

\noindent we can associate an undirected graph $G_L$ such that $\bb{L}$ is its Laplacian matrix.\footnote{Given a graph with adjacency matrix \bb{M}, its Laplacian matrix is defined as $\bb{L} = \bb{S} - \bb{M}$, where $\bb{S}$ is the degree matrix, i.e., a matrix having the degree of each node in the corresponding diagonal element, and zeros elsewhere.} Moreover, since at least $2k+1$ diagonal elements of $\bb{L}$ are non-zero (and strictly positive), at least $2k+1$ elements off the diagonal are non-zero (and strictly negative), since each row has to add up zero. 

Now, each off diagonal element of the laplacian matrix $\bb{L}$ corresponds to an edge of the graph $G_L$. Since the matrix $\bb{L}$ is symmetric, the graph $G_L$ is undirected, and each edge appears twice in the Laplacian matrix. Since there are at least $2k+1$ non-zero elements off the diagonal of $\bb{L}$, the auxiliary graph $G_L$ has at least $\left \lfloor \dfrac{2k+1}{2}\right \rfloor= k+1$ edges. It is easy to see that, if the number of edges is $e\geq k+1$, then the auxiliary graph $G_L$ has at most $C \leq n-(k+1)$ connected components. 

Remembering that the number of connected components $C$ is given by the multiplicity of the $0$ eigenvalue of the Laplacian matrix $\bb{L}$, then the rank of $\bb{L}$ has to be at least $rank(\bb{L}) = n- C \geq k+1$.

However, by construction, $\bb{L}$ is the sum of $k$ rank-one matrices, and therefore $rank(\bb{L})\leq k$, which is a contradiction. This proves that the only solution to problem \eqref{problemF} is the identity, concluding the proof of the theorem. $\blacksquare$
%\end{proof}

\vspace{0.4cm}

This proof, and therefore the class of equivalence between problems \eqref{eq:GM} and \eqref{eq:RGM}, can be further extended by noting that we originally asked for every eigenvector orthogonal to $\ones$ to have at least $2k+1$ non-zero elements, but in the proof we only used this fact for one of these eigenvectors. It might be the case, for instance, that only one of the elements $\bb{F}_{i,i}$ is different from one, and therefore the matrix $\bb{L}$ has rank one. In this case, it would be sufficient to ask for the corresponding eigenvector to have at least $3$ non-zero elements. The problem is that we do not know in advance which or how many of the eigenvectors will be used in the summation to construct the $\bb{L}$ matrix. However, it is possible to weaken the hypothesis as shown next.

Let us consider all the $k$ eigenvectors satisfying $u_i^T\ones =0$, and sort them according to the number of non-zero elements, such that $|u_1|_0 \leq |u_2|_0 \leq \dots \leq |u_k|_0$, where $|\cdot|_0$ is the $\ell_0$ pseudo-norm, which counts the non-zero elements of a vector.

\begin{theorem}
%\begin{thm}
\label{thik}
If $\bb{A}$ has no repeated eigenvalues, and there are $k$ eigenvectors $u_i$ such that $u_i^T\ones =0$, sorted as above and satisfying $|u_i|_0\geq 2i+1$,
then problems \eqref{eq:GM} and \eqref{eq:RGM} are equivalent.
%\end{thm}
\end{theorem}

%\noindent \textbf{Proof of Theorem \ref{thik}}
\noindent \textit{Proof:}

%\begin{proof}

This proof follows exactly the same procedure as the previous one, with minor changes.

Let us assume, as before, that there is a solution $\bb{F}\neq \bb{I}$ to problem \eqref{problemF}. In order to fulfill condition c1), the last $n-k$ diagonal elements of $F$ have to be one, i.e., $\bb{F}_{i,i}=1$ for $i=k+1,\ldots n$. For the first $k$ diagonal elements, there might be some $0$ values. Let $M$ be the greatest index of the eigenvectors actually used in the sum, meaning $M = \max \{i \in 1\ldots k : \bb{F}_{i,i} \neq 0\}$. 

We can then write 
$$\bb{L}=\sum_{i=1}^{M}(1-\bb{F}_{i,i})u_iu_i^T .$$

Since $|u_M|_0 \geq 2M+1$, the auxiliary graph $G_L$ has at least $M+1$ edges. Therefore the number of connected components satisfies $C \leq n-(M+1)$, and hence $rank(\bb{L}) \geq M+1$. The contradiction, as before, comes from the fact that $\bb{L}$ is the sum of $M$ rank one matrices, concluding the proof. $\blacksquare$
%\end{proof}

\vspace{0.4cm}

As noted above, a necessary condition for the problems \eqref{eq:GM} and \eqref{eq:RGM} to  be equivalent is for the automorphism group of $G_A$ to be the trivial group. Therefore, we have the following corollary:

\begin{corollary}
%\begin{cor}
\label{cor1}
%If $A$ has no repeated eigenvalues, and there are $k$ eigenvectors $u_i$ such that $u_i^T\ones =0$, sorted as above and satisfying $|u_i|_0\geq 2i+1$,
If $\bb{A}$ has no repeated eigenvalues, and there are $k$ eigenvectors $u_i$ such that $u_i^T\ones =0$, sorted according to their $\ell_0$ norm as above, and satisfying $|u_i|_0\geq 2i+1$,
then the automorphism group of the corresponding graph $G_A$ is the trivial group.
%\end{cor}
\end{corollary}

%\noindent \textbf{Proof of Corollary \ref{cor1}}
%
%$\blacksquare$

\section{Interpretation and additional results}
\label{additional}

It is clear that the spectral decomposition of an adjacency matrix provides a lot of information about the automorphism group of a graph, and the graph matching problem itself. However, very little is known about how the eigenvalues and eigenvectors of the adjacency matrix affect the graph properties. In this section, we discuss some links between these two fields, paying particular attention to the equivalence of problems \eqref{eq:GM} and \eqref{eq:RGM}, and also discussing more general novel properties.%, which are new up to our knowledge.

As noted in the previous section, asymmetry of a graph is a necessary condition for problems \eqref{eq:GM} and \eqref{eq:RGM} to be equivalent, although is not sufficient, with asymmetric regular graphs serving as counter-examples (see red region in Figure \ref{fig:bubble}). 
%\red{As observed above, a necessary condition for problems \eqref{eq:GM} and \eqref{eq:RGM} to be equivalent, is for the automorphism group of the graph to be the trivial group. However, this condition is not sufficient. Take for instance a regular graph, i.e., a graph such that each vertex has the same number of neighbors, and denote by $\bb{J}$ the barycenter of the set of doubly stochastic matrices, $\bb{J}=\dfrac{1}{n}\ones\ones^T$. Hence, it is very easy to see that, if $\bb{A}$ is the adjacency matrix of the regular graph, then $\bb{AJ}=\bb{JA}$. Therefore, there is a solution to problem \eqref{eq:RGM} which is not a permutation matrix. Since there are regular graphs with trivial automorphism group (like the Frucht graph for instance), then this condition cannot be sufficient.} \red{REPETIDO}
It is interesting to note that regular graphs have the vector $\ones$ as an eigenvector, and since the adjacency matrix is symmetric, its eigenvectors are orthogonal to each other, therefore there are $n-1$ eigenvectors satisfying $u_i^T\ones = 0$. Hence, the condition asked in \cite{alex} is violated not only by one eigenvector, but by $n-1$ of them. 

Besides this observation, it is not clear the interpretation of the non existence of eingenvectors perpendicular to $\ones$.

Let us focus now on the properties of eigenvectors orthogonal to $\ones$ with restricted support, as in the statement of Theorem \ref{th2k}.

%\noindent \textbf{The simplest case: one eigenvector $u$ such that $u^T\ones=0$}
\subsection{The simplest case: one single eigenvector $u$ such that $u^T\ones=0$}
\label{caso1}

Let us assume here that $\bb{A}$ is the adjacency matrix of a graph $G_A$ with no repeated eigenvalues and only one eigevector $u$ satisfying $u^T\ones=0$. Now, if this vector $u$ has strictly more than two non-zero entries, i.e., $|u|_0 > 2$, then this graph falls into the hypothesis of Theorem \ref{th2k}. Therefore, the graph has trivial automorphism group, and if $G_B$ is an isomorphic graph, problems \eqref{eq:GM} and \eqref{eq:RGM} are equivalent.

Since the sum of the entries of $u$ is zero, the only remaining case is when $u$ has exactly two non-zero elements. Assuming that the eigenvectors are normalized, the eigenvector $u$ is of the form 

$$u=\left(0,\ldots,0,\dfrac{1}{\sqrt{2}},0,\ldots,0,\dfrac{-1}{\sqrt{2}},0,\ldots,0\right) .$$

Let $s$ and $t$ be the indices of the non-zero coefficients. Since $u$ is an eigenvalue, then we have $\bb{A}u=\lambda u$. Now, denoting by $\bb{A}_s$ and $\bb{A}_t$ the columns of $\bb{A}$ at positions $s$ and $t$ respectively, and taking into account the particular structure of $u$, the product $\bb{A}u$ is simply the difference between these two columns: $\bb{A}u = \dfrac{1}{\sqrt{2}}(\bb{A}_s-\bb{A}_t) = \lambda u$. Therefore, columns $\bb{A}_s$ and $\bb{A}_t$ are identical, except for the coordinates $s$ and $t$. This means that the nodes corresponding to indices $s$ and $t$ have exactly the same connectivity pattern with the rest of the nodes in the graph. 

%Consider now the rest of the involved entries, which are $(s,s)$, $(s,t)$, $(t,s)$ and $(t,t)$ of matrix $\bb{A}$. It is easy to see that, since the entries in $\bb{A}$ are either $1$ or $0$, 
Consider now the rest of the involved entries (nodes). Let $\bb{\tilde{A}}$ be the $2\times 2$ sub-matrix formed by entries $(s,s)$, $(s,t)$, $(t,s)$ and $(t,t)$ of matrix $\bb{A}$, and let $w=\left(\frac{1}{\sqrt{2}},\frac{-1}{\sqrt{2}}\right)$. We have then $\bb{\tilde{A}}w=\lambda w$. It is easy to see that, since the entries in $\bb{\tilde{A}}$ are either $1$ or $0$, then only three values of $\lambda$ are possible: $-1$, $0$ and $1$, corresponding to the following situations: 

\begin{itemize}
\item $\lambda = -1$: the matrix is $\tilde{\bb{A}}= \begin{pmatrix}
0 & 1\\
1 & 0
\end{pmatrix}$, therefore nodes $s$ and $t$ are connected and they have no loops;
\item $\lambda = 0$: the matrix is either $\tilde{\bb{A}}= \begin{pmatrix}
1 & 1\\
1 & 1
\end{pmatrix}$ or 
$\tilde{\bb{A}}= \begin{pmatrix}
0 & 0\\
0 & 0
\end{pmatrix}$. Therefore, nodes $s$ and $t$ are either connected and both have loops, or not connected without loops;
\item $\lambda = 1$: the matrix is $\tilde{\bb{A}}= \begin{pmatrix}
1 & 0\\
0 & 1
\end{pmatrix}$, therefore nodes $s$ and $t$ are not connected and both of them have loops.
\end{itemize} 
 
Taking into account that nodes $s$ and $t$ have the same connectivity pattern with the rest of the graph, in any of the situations listed above, nodes $s$ and $t$ are interchangeable, meaning that there exists a non trivial automorphism of the graph $G_A$, namely, the automorphism which permutes nodes $s$ and $t$, and leaves the rest of the nodes unchanged. Therefore, for graphs with the corresponding adjacency matrix having a single eigenvector orthogonal to the unity vector, and this eigenvector having exactly two non-zero entries, problems \eqref{eq:GM} and \eqref{eq:RGM} are note equivalent. The problems are equivalent if the eigenvector has more than two non-zero entries. 
%coordinates $s$ and $t$ of columns $\bb{A}_s$ and $\bb{A}_t$. 

%\noindent \textbf{The general case}
\subsection{The general case}

Let us further analyze the relationship between the group of automorphisms of a graph and the eigenvectors of its adjacency matrix, now considering a more general case in terms of the non-zero elements of the eigenvectors.

First, observe that if the matrix $\bb{A}$ has simple spectrum, then each element of the automorphism group has order two, with the exception of the identity (this result appears in \cite{biggs} and \cite{lovasz}): %. An alternative way to state this is like follows:

\begin{lemma}
If $\bb{A}$ has no repeated eigenvalues and $\bb{P}$ is a permutation matrix such that\\ $\bb{AP}=\bb{PA}$, then $\bb{P}^2=\bb{I}$.
\end{lemma}
%\noindent \textbf{Proof}
\noindent \textit{Proof:}
%\begin{proof}

In order to prove this, let $u$ be an eigenvector of $\bb{A}$ associated with the eigenvalue $\lambda$. Then, $\bb{A}\bb{P}u = \bb{PA}u = \bb{P}\lambda u = \lambda \bb{P}u$. Therefore, the vector $\bb{P}u$ is an eigenvector associated with the eigenvalue $\lambda$ as well. Since every eigenspace has dimension $1$, and the multiplication by the permutation matrix preserves the norm, then necessarily $\bb{P}u = \pm u$, and hence $\bb{P}^2u=u$. Since this is true for every eigenvector $u$ in the basis, then $\bb{P}^2=\bb{I}$.  $\blacksquare$
%\end{proof}

We are now able to prove the following.

%\begin{prop} If $\bb{A}$ has no repeated eigenvalues, and the group of automorphisms of $G_A$ is not trivial, 
%then there exist a set of $k$ eigenvectors $u_i$ satisfying $u_i^T\ones =0$, each one of them with exactly $2k$ non-zero entries, in the same location.
%\label{prop1}
%\end{prop}

\begin{proposition}
%\begin{prop}
 If $\bb{A}$ has no repeated eigenvalues, and the group of automorphisms of $G_A$ is non trivial, 
then there exist a set of $k$ eigenvectors $u_i$ satisfying $u_i^T\ones =0$, each one of them having at most $2k$ non-zero entries.
\label{prop1}
%\end{prop}
\end{proposition}

%\noindent \textbf{Proof of Proposition \ref{prop1}}
\noindent \textit{Proof:}
%\begin{proof}

Let $\bb{P} \neq \bb{I}$ be a permutation matrix, corresponding to a non-trivial automorphism of $G_A$. As observed above, since $\bb{A}$ has simple spectrum, then $\bb{P}^2=\bb{I}$. Since the permutation has order two, we can re-arrange the order of the nodes in such a way that the resulting permutation matrix $\bb{P}$ is block diagonal as follows %in Figure \ref{perm_block}.

$$
 \bb{P}= \left(
    \begin{array}{cccccccc}
0 & 1 & &   &    &     &    & \text{\huge{0}}   \\ %\cline{1-1}
1 & 0 & &   &    &     &    &  \\ %\cline{2-2}
& &  \ddots &           &  &  &      &  \\ %\cline{4-4}
 &  & &  0& 1  &    &        &  \\ %\cline{2-2}
&  &  &  1& 0  &    &         &  \\ %\cline{2-2}
& & &           &    & 1    &     &  \\ %\cline{3-3}
& & &           &  &  &  \ddots    &  \\ %\cline{4-4}
\text{\huge{0}} & & &           &          &       &  & 1 \\ %\cline{5-5}
  \end{array}\right).
$$

As in the previous section, consider the eigen-decomposition $\bb{A} = \bb{UDU^T}$, which transforms the problem $\bb{AP}=\bb{PA}$ into $\bb{DF}=\bb{FD}$, where the new unknown matrix $\bb{F}$ is defined as $\bb{F} = \bb{U^TPU}$, or equivalently, $\bb{P}=\bb{UFU^T}$. As before, since $\bb{A}$ has no repeated eigenvalues, $\bb{F}$ is necessarily diagonal, and therefore $\bb{P}=\bb{UFU^T}$ is one possible eigen-decomposition of $\bb{P}$.

Now, the matrix $\bb{U}$ of normalized eigenvectors of $\bb{A}$ is unique, up to changes of sign in each column. This is not true for $\bb{P}$, since it has repeated eigenvalues. However, any orthogonal eigen-decomposition of $\bb{P}$ can be obtained as an orthogonal transformation (rotation and/or symmetries) of $\bb{U}$.

Given that $\bb{P}$ is block-diagonal, one possible eigen-decomposition can be obtained by combining the eigen-decompositions of each block. The lower part of $\bb{P}$ is an identity block, and hence all eigenvalues are equal to $1$, with canonical eigenvectors. The rest are $2 \times 2$ blocks with the following decomposition:

$$
\left(
 \begin{array}{cc}
0 & 1 \\
1 & 0 \\    
\end{array}\right) = 
\left(
 \begin{array}{cc}
\frac{-1}{\sqrt{2}} & \frac{1}{\sqrt{2}} \\
\frac{1}{\sqrt{2}} & \frac{1}{\sqrt{2}} \\    
\end{array}\right)
\left(
 \begin{array}{cc}
-1 & 0 \\
0 & 1 \\    
\end{array}\right)
\left(
 \begin{array}{cc}
\frac{-1}{\sqrt{2}} & \frac{1}{\sqrt{2}} \\
\frac{1}{\sqrt{2}} & \frac{1}{\sqrt{2}} \\    
\end{array}\right) .
$$

Therefore, a plausible eigendecomposition for $\bb{P}$ is $\bb{P} = \bb{VEV^T}$, with:

$$
 \bb{E}= \left(
    \begin{array}{cccccccc}
-1 &  & &   &    &     &    & \text{\huge{0}}   \\ %\cline{1-1}
 & 1 & &   &    &     &    &  \\ %\cline{2-2}
& &  \ddots &           &  &  &      &  \\ %\cline{4-4}
 &  & &  -1&   &    &        &  \\ %\cline{2-2}
&  &  &  & 1  &    &         &  \\ %\cline{2-2}
& & &           &    & 1    &     &  \\ %\cline{3-3}
& & &           &  &  &  \ddots    &  \\ %\cline{4-4}
\text{\huge{0}} & & &           &          &       &  & 1 \\ %\cline{5-5}
  \end{array}\right)
 \quad  , \quad
 \bb{V}= \left(
    \begin{array}{cccccccc}
\frac{-1}{\sqrt{2}} & \frac{1}{\sqrt{2}} & &   &    &     &    & \text{\huge{0}}   \\ %\cline{1-1}
\frac{1}{\sqrt{2}} & \frac{1}{\sqrt{2}} & &   &    &     &    &  \\ %\cline{2-2}
& &  \ddots &           &  &  &      &  \\ %\cline{4-4}
 &  & &  \frac{-1}{\sqrt{2}} & \frac{1}{\sqrt{2}}  &    &        &  \\ %\cline{2-2}
&  &  &  \frac{1}{\sqrt{2}} & \frac{1}{\sqrt{2}}  &    &         &  \\ %\cline{2-2}
& & &           &    & 1    &     &  \\ %\cline{3-3}
& & &           &  &  &  \ddots    &  \\ %\cline{4-4}
\text{\huge{0}} & & &           &          &       &  & 1 \\ %\cline{5-5}
  \end{array}\right) .
$$

%Hence, the matrix $\bb{U}$ can be written as $\bb{U}=\bb{TV}$, where $\bb{T}$ is an orthogonal matrix that leaves invariant the eigenspaces $S_{-1}$ and $S_{1}$.
Since the columns of both $\bb{U}$ and $\bb{V}$ are possible basis of the eigenvalues of $\bb{P}$, the matrix $\bb{U}$ can be thought as an orthogonal transformation of $\bb{V}$ that leaves invariant the eigenspaces $S_{-1}$ and $S_{1}$ (eigenspaces of $\bb{P}$ associated with eigenvalues $-1$ and $1$ respectively).
 Observe that the eigenspace $S_{-1}$ is composed by vectors orthogonal to $\ones$, and therefore, in the orthogonal transformation from $\bb{V}$ to $\bb{U}$, the whole subspace $S_{-1}$ will be mapped to a subspace orthogonal to $\ones$. 
Let $k$ be the dimension of the subspace $S_{-1}$, and let us denote by $\tilde{\bb{U}}$ the set of the $k$ eigenvectors of $\bb{A}$ (columns of $\bb{U}$) corresponding to the eigenspace $S_{-1}$ after the linear mapping. These columns of $\bb{U}$, as argued above, are orthogonal to $\ones$. Analogously, let $\tilde{\bb{V}}$ be formed by the columns of $\bb{V}$ associated with the eigenvalue $-1$, so the columns of $\tilde{\bb{V}}$ are a basis of $S_{-1}$.

Given that we assumed $\bb{P}\neq \bb{I}$, there is at least one $2\times 2$ non identity block like the one described above, and therefore $k \geq 1$. Since $S_{-1}$ is invariant under the orthogonal transformation, then the mapping of this subspace can be written as linear combinations of the elements of the basis, this is, $\tilde{\bb{U}} = \tilde{\bb{V}}\bb{T}$, where $\bb{T}$ is an orthogonal matrix. %, $\tilde{\bb{V}}$ is a basis of $S_{-1}$, and $\tilde{\bb{U}}$ is the set of $k$ eigenvectors of $A$, orthogonal to $\ones$. 
Since, according to the previous description, $\tilde{\bb{V}}$ has the form

$$ \tilde{\bb{V}} = \frac{1}{\sqrt{2}} 
\underbrace{
\left(
    \begin{array}{rrrr}
    -1 & 0 & \ldots & 0\\
     1 & 0 & \ldots & 0\\
    0 & -1 & \ldots & 0\\
     0 & 1 & \ldots & 0\\
     \vdots & \vdots &  & \vdots\\        
     0 & 0 & \ldots & -1\\     
     0 & 0 & \ldots & 1\\     
          \vdots & \vdots &  & \vdots\\        
%	0 & 0 & \ldots & 0\\               
%    \underbrace{0 & 0 & \ldots & 0}_{1} \\     
  \end{array}\right)}_{k \text{ columns}} ,    
$$

then $\tilde{\bb{U}}$ is conformed by $k$ vectors with $2k$ non-zero entries at most. Moreover, each one of these vectors has an even number of non-zero entries.
 $\blacksquare$
 %\end{proof}
 
 \vspace{0.4cm}

We have simulated millions of Erd\H{o}s-R\'enyi graphs, and obtained have empirical evidence suggesting that there is always a subset of these $k$ vectors formed by $r$ vectors with exactly $2r$ non-zero entries each, in the same location,  which correspond to the $2r$ nodes which are permuted in the automorphism (see the examples in the appendix). However, the arguments used in the previous proof are not sufficient to formally prove this. 
 
 On the other hand, the empirical evidence also suggests that a converse of this last statement may be true. We formulate then the following conjecture.
 
 \begin{conjecture}
%\red{
%\begin{conj} 
If $\bb{A}$ has no repeated eigenvalues, and there exist a set of $r$ eigenvectors $u_i$ satisfying $u_i^T\ones =0$, each one of them with exactly $2r$ non-zero entries, in the same location, 
then the group of automorphisms of $G_A$ is not trivial.
\label{conj1}
%\end{conj}
 \end{conjecture}
%}

We know that this is true for $k=1$, %since the proof for this particular case is contained in Section \ref{caso1}. 
%since this is what has been shown in the simplest example of 
Section \ref{caso1}. 
The proof for the general case, as well as other relations between spectral properties and the automorphism group,  are part of future work.

\section{Graph Examples}
\label{graph_examples}
Figure \ref{fig:bubble} shows different sets of graphs according to the relevant characteristics for this paper, principally about eigenvectors and eigenvalues. For instance, the class of graphs where theorems \ref{th2k} and \ref{thik} apply lays on the intersection of asymmetric and simple spectrum graphs, but outside the \textit{non-orthogonal to $\ones$} set. It is important therefore to show that there exist graphs in this subset, and in general that each subset in the diagram is not empty.

\begin{figure}[h!]
\centering
\hspace{2cm} \parbox{0.4\textwidth}{
 \begin{tikzpicture}[thick,scale=0.7]
    \node (1) at (    0.8244 ,   2.5802 )  {};
    \node (2) at (    1.3403 ,   1.7597 )  {};
    \node (3) at (    1.8235 ,   2.7188 )  {};
    \node (4) at (    3.2388 ,   4.1572 )  {};
    \node (5) at (    5.8952 ,   4.0242 )  {};
    \node (6) at (    7.0983 ,   2.8054 )  {};
    \node (7) at (    7.4529 ,   1.7837 )  {};
    \node (8) at (    8.0022 ,   2.6235 )  {};
    \node (9) at (    6.9794 ,   5.3356 )  {};
    \node (10) at (    7.9368 ,   5.5782 )  {};
    \node (11) at (    7.2887 ,   6.3242 )  {};
    \node (12) at (    1.8911 ,   6.1421 )  {};

\draw (1) to[out=100,in=225,looseness=1] (12);
\draw (1) to[out=270,in=150,looseness=0.6] (2);
\draw (1) to[out=50,in=150,looseness=0.6] (3);
\draw (2) to[out=30,in=280,looseness=0.6] (3);
\draw (2) to[out=310,in=230,looseness=0.6] (7);
\draw (7) to[out=20,in=290,looseness=0.6] (8);
\draw (6) to[out=240,in=140,looseness=0.6] (7);
\draw (6) to[out=40,in=120,looseness=0.6] (8);
\draw (8) to[out=40,in=310,looseness=0.6] (10);
\draw (10) to[out=90,in=340,looseness=0.6] (11);
\draw (10) to[out=230,in=350,looseness=0.6] (9);
\draw (11) to[out=200,in=120,looseness=0.6] (9);
\draw (11) to[out=140,in=40,looseness=0.6] (12);
%\draw (3) -- (4) -- (5) -- (6);
\draw (3) to[out=20,in=270,looseness=0.6] (4);
\draw (4) to[out=340,in=190,looseness=0.6] (5);
\draw (5) to[out=280,in=160,looseness=0.6] (6);
%\draw (5) -- (9);
\draw (5) to[out=95,in=190,looseness=0.6] (9);
%\draw (4) -- (12);
\draw (4) to[out=100,in=340,looseness=0.6] (12);
\end{tikzpicture}
}
\quad
\begin{minipage}{0.4\textwidth}
\begin{tikzpicture}[thick,scale=0.5]
        \node (1) at (	  8.5495  ,  2.1965 ) {};
        \node (2) at (    1.4505  ,  7.8035 ) {};
        \node (3) at (    6.3557  ,  9.5362 ) {};
        \node (4) at (    9.3874  ,  5.0000 ) {};
        \node (5) at (    8.5495  ,  7.8035 ) {};
        \node (6) at (    6.3557  ,  0.4638 ) {};
        \node (7) at (    3.6443  ,  0.4638 ) {};
        \node (8) at (    0.6126  ,  5.0000 ) {};
        \node (9) at (    3.6443  ,  9.5362 ) {};
        \node (10) at (    1.4505  ,  2.1965 ) {};                

 \draw (1) -- (4) -- (5) -- (3)  -- (9)  -- (2) -- (8) -- (10) -- (7) -- (6) -- (1);
 \draw (1) -- (5) -- (8) -- (7) -- (3) -- (6) -- (2) -- (10) -- (4) --(9) -- (1);
\end{tikzpicture}

\end{minipage}

\caption{Regular graphs with trivial automorphism group. Left: the Frucht graph, with simple spectrum. Right: a regular asymmetric graph with repeated eigenvalues.} \label{fig:frucht}
\end{figure}
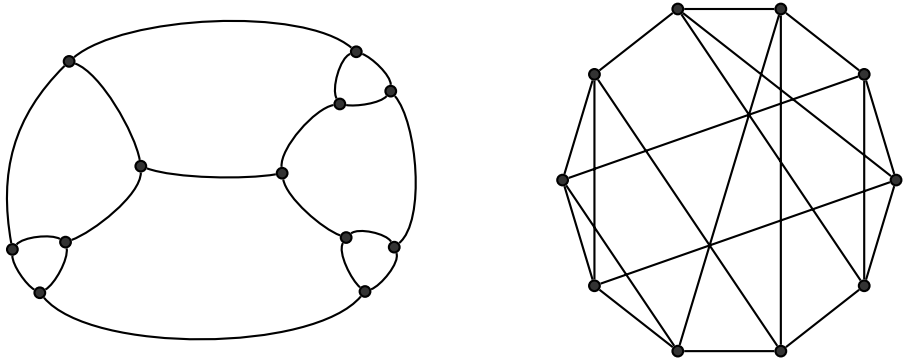

As mentioned above, the Frucht graph \cite{frucht}, illustrated in Figure \ref{fig:frucht} (left), serves as an example of regular graphs with trivial automorphism group and simple spectrum. The regular graph in Figure \ref{fig:frucht} (right) also has trivial automorphism group, but the adjacency matrix has repeated eigenvalues.

\tikzstyle{every node}=[circle, draw, fill=black!80, inner sep=0pt, minimum width=4pt]

%\begin{figure}[h!]
%\centering
%\caption{regular 10 nodos}\label{fig:graph_zone_8}
%\end{figure}

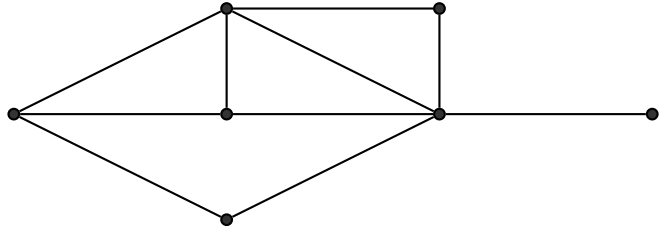
\begin{figure}[h!]
\centering

\begin{tikzpicture}[thick,scale=0.7]
    \node (1) at (8,0) {};
    \node (2) at (4, 0)  {};
    \node (3) at (0, -2) {};
    \node (4) at (-4, 0) {};
    \node (5) at (0, 2) {};
    \node (6) at (4, 2) {};
    \node (7) at (0, 0) {};

 \draw (1) -- (2) -- (3) -- (4)  -- (5)  -- (6) -- (2) -- (5) -- (7) -- (4);
 \draw (7) -- (2);
\end{tikzpicture}
\caption{A non-friendly graph, with simple spectrum but one eigenvector orthogonal to $\ones$.}\label{fig:graph_zone_3}
\end{figure}

The graph of Figure \ref{fig:graph_zone_3} is an example of a graph with simple spectrum, but where there is an eigenvector $u$ such that $u^T\ones=0$, and therefore this is not a \textit{friendly} graph. This eigenvector has $4$ non-zero elements, and hence Theorem \ref{th2k} applies. As a consequence, for any isomorphic graph, problems \eqref{eq:GM} and \eqref{eq:RGM} are equivalent, and in particular the automorphism group of this graph is trivial. Since the graph is not \textit{friendly}, the results of \cite{alex} do not hold for this graph.

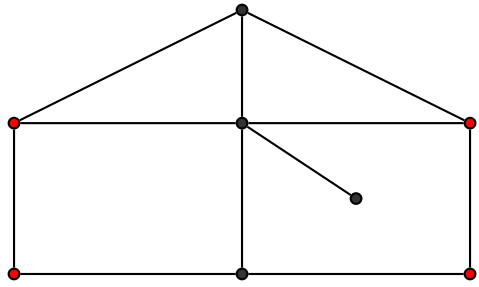
\begin{figure}[h!]
\centering
\begin{tikzpicture}[thick,scale=1]
    \node (1) at (0,1.5) {};
    \node[fill=red] (2) at (3, 0)  {};
    \node[fill=red] (3) at (3, -2) {};
    \node (4) at (0, -2) {};
    \node[fill=red] (5) at (-3, -2) {};
    \node[fill=red] (6) at (-3, 0) {};
    \node (7) at (0, 0) {};
    \node (8) at (1.5, -1) {};

 \draw (1) -- (2) -- (3) -- (4)  -- (5)  -- (6) -- (1) -- (7) -- (8);
 \draw (6) -- (7) -- (2);
 \draw (4) -- (7);
\end{tikzpicture}
\caption{A graph with simple spectrum but non-trivial automorphism group.}\label{fig:graph_zone_4}
\end{figure}

The graph in Figure \ref{fig:graph_zone_4} has simple spectrum but non-trivial automorphism group. Indeed, it has two eigenvectors $u_1$ and $u_2$ orthogonal to $\ones$, each one with four non-zero elements. Namely, the eigenvectors are 
$$u_1 = \left(\frac{1}{2},\frac{1}{2},\frac{-1}{2},\frac{-1}{2},0,0,0,0\right)$$
and
$$u_2 = \left(\frac{1}{2},\frac{-1}{2},\frac{1}{2},\frac{-1}{2},0,0,0,0\right) .$$

Illustrating the Conjecture \ref{conj1}, the first four coordinates of the eigenvectors correspond to the four red nodes of the graph. The non-trivial automorphism consist on permuting the two lower nodes between themselves, and the two upper nodes between themselves, as it can be clearly seen in the figure. Of course, theorems \ref{th2k} and \ref{thik} do not apply for this graph.

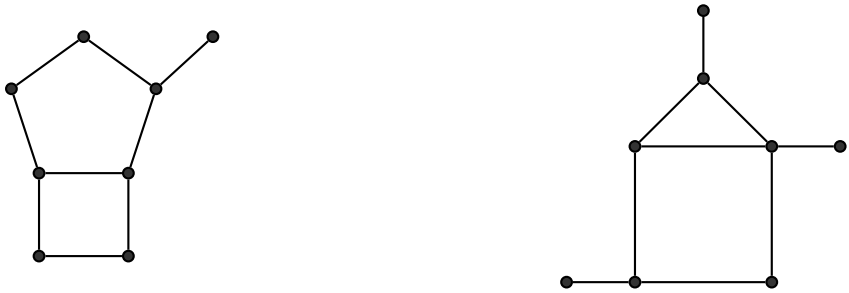
\begin{figure}[h!]
%\centering
\hspace{3cm}
\parbox{0.4\textwidth}{
\begin{tikzpicture}[thick,scale=1]
    \node (1) at (1.7,1) {};
    \node (2) at (0.9511,0.3090)  {};
    \node (3) at (0, 1) {};
    \node (4) at (-0.9511,0.3090) {};
    \node (5) at (-0.5878,-0.8090 ) {};
    \node (6) at (-0.5878,-1.9090 ) {};
    \node (7) at (0.5878,-1.9090 ) {};
    \node (8) at (0.5878,-0.8090 ) {};

 \draw (1) -- (2) -- (3) -- (4)  -- (5) -- (6) -- (7) -- (8) -- (2);
 \draw (5) -- (8);
\end{tikzpicture}

}
\quad
\begin{minipage}{0.4\textwidth}
\begin{tikzpicture}[thick,scale=0.9]
    \node (1) at (-2,-2) {};
    \node (2) at (-1, -2)  {};
    \node (3) at (1, -2) {};
    \node (4) at (1, 0) {};
    \node (5) at (2, 0) {};
    \node (6) at (0, 1) {};
    \node (7) at (0, 2) {};
    \node (8) at (-1, 0) {};

 \draw (1) -- (2) -- (3) -- (4)  -- (5);
 \draw (2) -- (8) -- (4) -- (6) -- (7);
 \draw (6) -- (8);
\end{tikzpicture}

\end{minipage}
\caption{Asymmetric graphs with repeated eigenvalues, with (left) and without (right) eigenvectors orthogonal to $\ones$.}\label{fig:graph_zones_5_1}
%\caption{An asymmetric graph with repeated eigenvalues and one eigenvector orthogonal to $\ones$.}\label{fig:graph_zone_1}
\end{figure}

The following two graphs, illustrated in Figure \ref{fig:graph_zones_5_1}, have trivial automorphism group and both have repeated eigenvalues. The first one has one eigenvector orthogonal to $\ones$, while the second one has no eigenvector $u$ satisfying $u^T\ones=0$.

\begin{figure}[h!]
\centering
\includegraphics[width=0.7\textwidth]{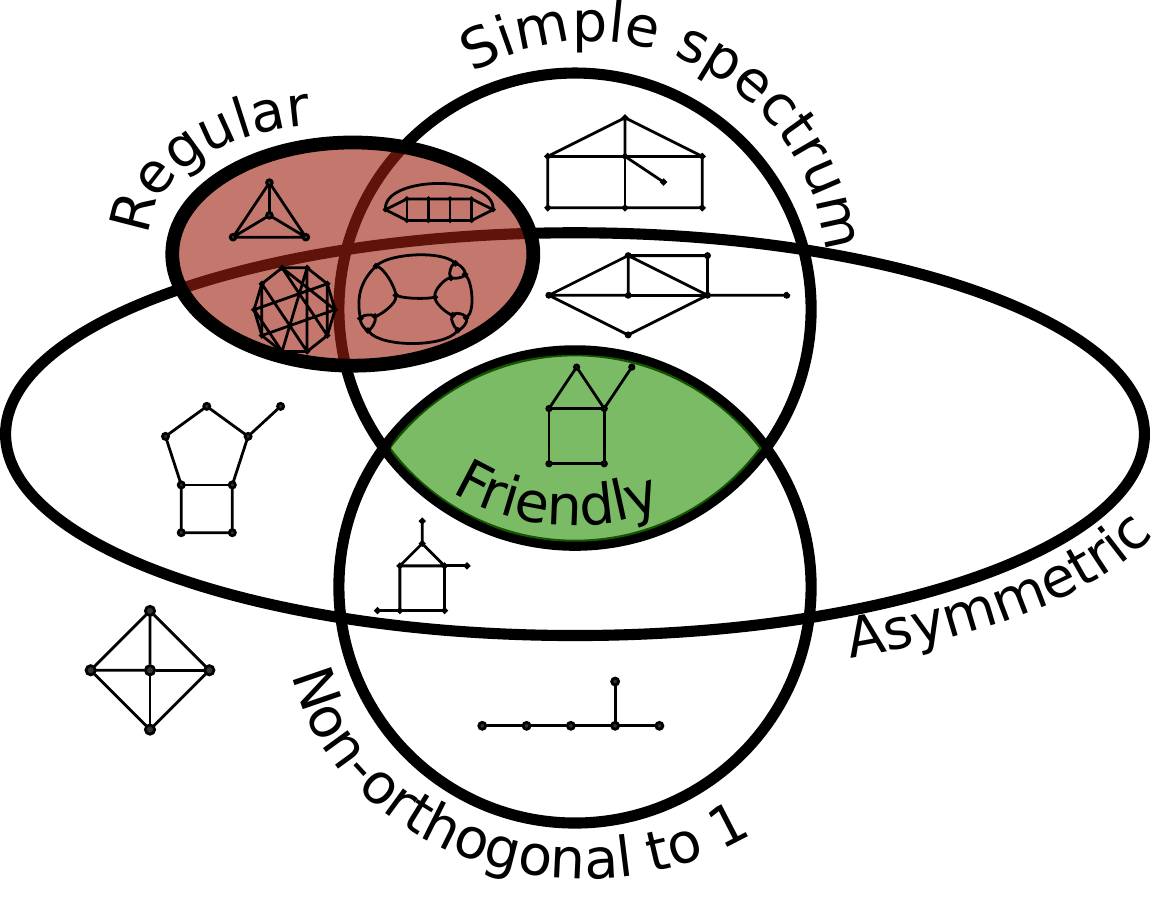} 
\caption{Examples of graphs in each class.}
\label{fig:bubble_con_grafos}
\end{figure}

Finally, Figure \ref{fig:bubble_con_grafos} shows the same diagram as in Figure \ref{fig:bubble}, but with examples of graphs inside each intersection, demonstrating that none of these subsets is empty. 

\section{Conclusion}

We have addressed the equivalence of the graph matching problem with its most common convex relaxation, generalizing the results in \cite{alex}, and extending the analysis to graph automorphism properties.

Theorem \ref{th2k} and the stronger version, Theorem \ref{thik}, state conditions on the spectral properties of the adjacency matrix of a graph in order for the graph matching problem and the convex relaxation to be equivalent. Specifically, if the adjacency matrix has simple spectrum, and the eigenvectors orthogonal to vector $\ones$ have enough non-zero entries, then the equivalence between the two problems holds. This gives also a set of easily verifiable conditions implying that the automorphism group of a graph is trivial.

\begin{figure}[h!]
\centering
\includegraphics[width=0.6\textwidth]{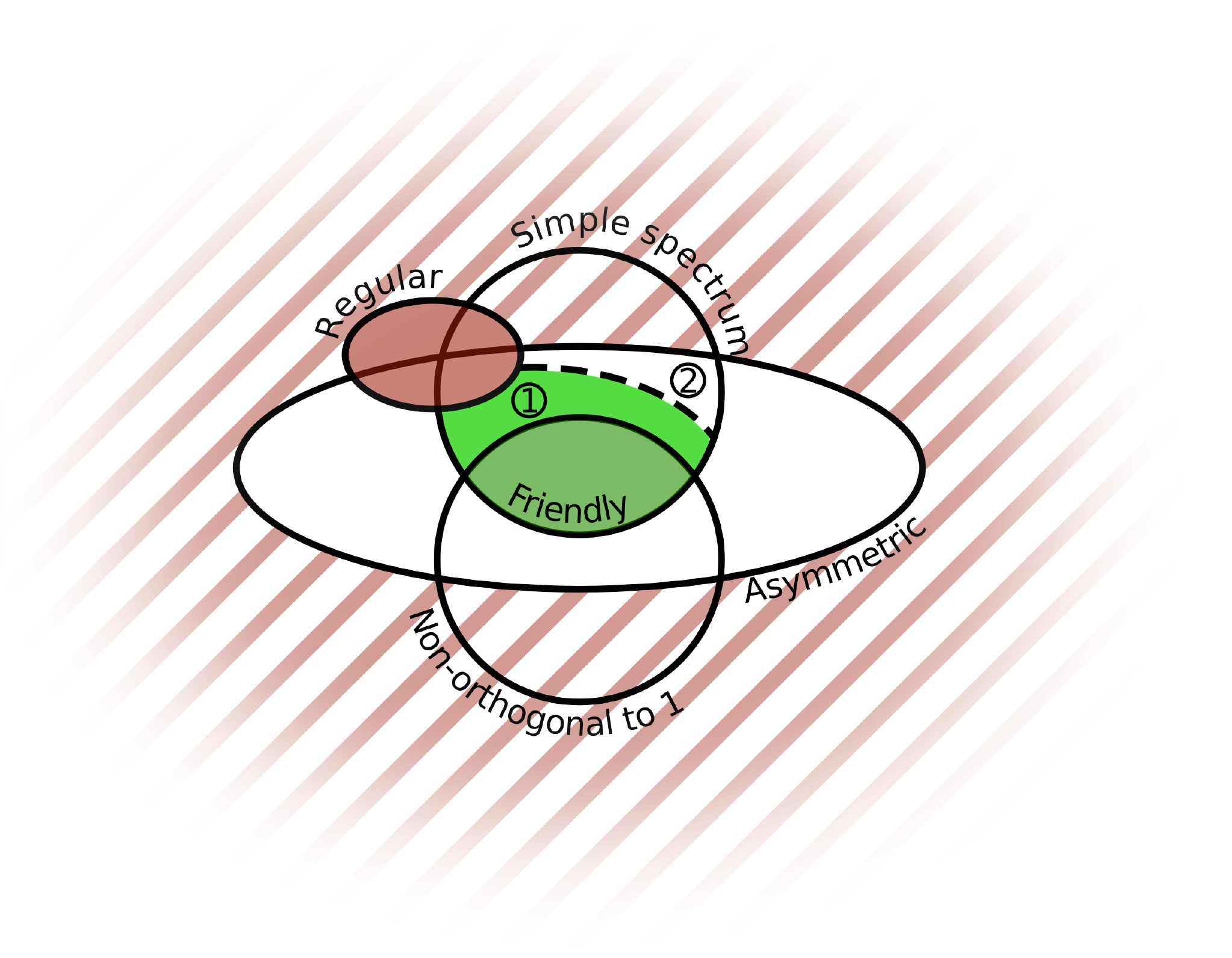} 
\caption{Regions where problems \eqref{eq:GM} and \eqref{eq:RGM} are known to be equivalent (green) and non-equivalent (red). Outside the \textit{asymmetric} set and inside the \textit{regular} set, the problems are known to be non-equivalent. Problems are equivalent for \textit{friendly} graphs \cite{alex}, and in the green zone $1$ by virtue of the theorems proved in Section \ref{sec:thms}. The zone $2$ consists of non-regular asymmetric graphs with simple spectrum, but not satisfying the conditions of theorems \ref{th2k} or \ref{thik}; although this subset might be nonempty, since we could not find examples of graphs in this zone.}
\label{fig:bubble_ext}
\end{figure}

The extension of the set where problems \eqref{eq:GM} and \eqref{eq:RGM} are known to be equivalent, due to these new results, is shown in Figure \ref{fig:bubble_ext}.

In addition to the main theorems, we provided evidence that these particular eigenvectors, orthogonal to $\ones$, contain critical information about the symmetries of the graph, specially in their non-zero entries.

During the last decades, important theory was developed on eigenvalues and eigenvectors of the Laplacian matrix of a graph, with very important theoretic results, which brought important applications. The new results here presented shed light on some spectral properties of the adjacency matrix, and leave open some other questions about the link between these properties and the automorphisms of the graph.

%\begin{figure}[h!]
%\centering
%\def\svgwidth{270pt}
%\begin{scriptsize}
%\input{bubble1.pdf_tex}
%\end{scriptsize}
%\caption{Mismatch .}
%\label{fig:energy}
%\end{figure}

%\begin{figure}[h!]
%\centering
%\includegraphics[width=0.49\textwidth]{bubble.pdf} 
%\caption{Mismatch.}
%\label{fig:bubble}
%\end{figure}

%\red{finish}%
%\noindent \textbf{Acknowledgments}

\section*{Acknowledgments}

M.F. performed part of this work while at Duke University.
We thank Dr. Pablo Sprechmann, Prof. Pablo Mus\'e, Prof. Alex Bronstein, and Prof. Josh Vogelstein for important feedback.
Work partially supported by NSF, ONR, NGA, ARO, AFOSR, and ANII.
%\newpage

\bibliographystyle{imaiai}
\bibliography{biblio}

\ifx\undefined\BySame
\newcommand{\BySame}{\leavevmode\rule[.5ex]{3em}{.5pt}\ }
\fi
\ifx\undefined\textsc
\newcommand{\textsc}[1]{{\sc #1}}
\newcommand{\emph}[1]{{\em #1\/}}
\let\tmpsmall\small
\renewcommand{\small}{\tmpsmall\sc}
\fi
\begin{thebibliography}{99}

\bibitem{alex}
\textsc{Aflalo, Y., Bronstein, A.  {\small \&} Kimmel, R.}  (2014) Graph
  matching: relax or not?. \emph{arXiv:1401.7623}.

\bibitem{biggs}
\textsc{Biggs, N.}  (1993) \emph{Algebraic Graph Theory}, Cambridge
  Mathematical Library. Cambridge University Press.

\bibitem{ConteReview}
\textsc{Conte, D., Foggia, P., Sansone, C.  {\small \&} Vento, M.}  (2004)
  Thirty years of graph matching in pattern recognition. \emph{International
  Journal of Pattern Recognition and Artificial Intelligence}, \textbf{18}(03),
  265--298.

\bibitem{fiori2013nips}
\textsc{Fiori, M., Sprechmann, P., Vogelstein, J., Mus\'{e}, P.  {\small \&}
  Sapiro, G.}  (2013) Robust Multimodal Graph Matching: Sparse Coding Meets
  Graph Matching. \emph{Advances in Neural Information Processing Systems 26},
  pp. 127--135.

\bibitem{frucht}
\textsc{Frucht, R.}  (1949) Graphs of degree three with a given abstract group.
  \emph{Canadian J. Math}, \textbf{1}, 365--378.

\bibitem{lovasz}
\textsc{Lov{\'a}sz, L.}  (2007) Eigenvalues of graphs. \emph{Lecture notes},
  \url{http://www.cs.elte.hu/~lovasz/eigenvals-x.pdf}.

\bibitem{vince}
\textsc{Lyzinski, V., Fishkind, D., Fiori, M., Vogelstein, J., Priebe, C.
  {\small \&} Sapiro, G.}  (2014) Graph Matching: Relax at Your Own Risk.
  \emph{arXiv preprint arXiv:1405.3133}.

\bibitem{umeyama}
\textsc{Umeyama, S.}  (1988) An eigendecomposition approach to weighted graph
  matching problems. \emph{Pattern Analysis and Machine Intelligence, IEEE
  Transactions on}, \textbf{10}(5), 695--703.

\bibitem{FAQ}
\textsc{Vogelstein, J., Conroy, J., Podrazik, L., Kratzer, S., Harley, E.,
  Fishkind, D., Vogelstein, R.  {\small \&} Priebe, C.}  (2012) Fast
  Approximate Quadratic Programming for Large (Brain) Graph Matching.
  \emph{arXiv:1112.5507}.

\bibitem{Zaslavskiy2009}
\textsc{Zaslavskiy, M., Bach, F.  {\small \&} Vert, J.}  (2009) A Path
  Following Algorithm for the Graph Matching Problem. \emph{Pattern Analysis
  and Machine Intelligence, IEEE Transactions on}, \textbf{31}(12), 2227--2242.

\end{thebibliography}

%\appendix

%\input{figurafinal.tex}

\end{document}